\documentclass{gtart}

\def\ifplaintex{\expandafter\ifx\csname documentclass\endcsname\relax}

\ifplaintex 
\else
\fi

\def\gtm{{\mathsurround=0pt\it $\cal G\mskip-2mu$eometry \&\ 
$\cal T\!\!$opology $\cal M\mskip-1mu$onographs}}    %  for monographs

\def\gtp{{\mathsurround=0pt\it $\cal G\mskip-2mu$eometry \&\ 
$\cal T\!\!$opology $\cal P\!$ublications}}  % GT publications

\def\recd{{\small Received:\qua\receiveddate\ifx\reviseddate\relax
\else\qquad Revised:\qua\reviseddate\fi\par}} 

\def\volumenumber#1{\def\thevolumenumber{#1}}
\def\volumeyear#1{\def\thevolumeyear{#1}}
\def\volumename#1{\def\thevolumename{#1}}
\def\papernumber#1{\def\thepapernumber{#1}}
\def\pagenumbers#1#2{\def\startpage{#1}\def\finishpage{#2}}
\def\published#1{\def\publishdate{#1}}
\def\received#1{\def\receiveddate{#1}}
\def\revised#1{\def\reviseddate{#1}}
\def\accepted#1{\def\accepteddate{#1}}
\def\asciititle#1{\def\theasciititle{#1}}

\long\def\asciiabstract#1{\long\def\theasciiabstract{#1}}

\def\shorttitle#1{\def\theshorttitle{#1}}

\let\\\par
\let\thevolumenumber\relax\let\thepapernumber\relax
\let\thevolumeyear\relax\let\startpage\relax
\let\finishpage\relax\let\publishdate\relax\let\receiveddate\relax
\let\reviseddate\relax\let\accepteddate\relax\let\theasciititle\relax
\let\theasciiauthors\relax
\let\theasciiabstract\relax

\let\theerratum\relax\let\theasciiemail\relax
\let\theshortauthors\relax\let\theshorttitle\relax

\def\startpage{1}\def\finishpage{15}\def\thepapernumber{77}

\volumenumber{2}
\volumename{Proceedings of the Kirbyfest}
\volumeyear{1999}

\long\def\maketitlep{   % start of definition of \maketitlep

\count0=\startpage

\gtm\nl        %   GT mongraphs (top left) 
{\small Volume \thevolumenumber: \thevolumename\nl 
\ifx\theerratum\relax\else Erratum \erratumnumber\nl\fi
Pages \startpage--\finishpage\nl}

\vglue 0.1truein   % top margin

{\parskip=0pt\leftskip 0pt plus 1fil\def\\{\par\smallskip}{\ifplaintex\large
\else\Large\fi\bf\thetitle}\par\medskip}   
\vglue 0.05truein 

{\parskip=0pt\leftskip 0pt plus 1fil\def\\{\par}{\sc\theauthors}
\par\medskip}%
 
\vglue 0.03truein 

{\small\leftskip 25pt\rightskip 25pt{\bf Abstract}\stdspace\theabstract

{\bf AMS Classification}\stdspace\theprimaryclass
\ifx\thesecondaryclass\relax\else; \thesecondaryclass\fi\par
{\bf Keywords}\stdspace \thekeywords\par}\vglue 7pt

}   % end of definition of \maketitlep

\font\phead=cmsl9 scaled 950
\font=cmsl9 scaled 1050
\font\pnum=cmbx10 scaled 913
\font\lnum=cmbx10 
\font\pfoot=cmsl9 scaled 950
\font=cmsl9 scaled 1050
\ifplaintex
\headline{\vbox to 0pt{\vskip -4.5mm\line{\small\phead\ifnum
\count0=\startpage ISSN 1464-8997 (on line)
1464-8989 (printed) \hfill {\pnum\folio}\else\ifodd\count0\def\\{ }% 
\ifx\theshorttitle\relax\thetitle\else\theshorttitle\fi\hfill{\pnum\folio}
\else\def\\{ and }{\pnum\folio}\hfill\ifx\theshortauthors\relax\theauthors
\else\theshortauthors\fi\fi\fi}\vss}}
\footline{\vbox to 0pt{\vglue 0mm\line{\small\pfoot\ifnum\count0=\startpage
Published \publishdate:\qua\copyright\ \gtp\hfill\else
\gtm, Volume \thevolumenumber\ (\thevolumeyear)\hfill\fi}\vss
}}
\else
\makeatletter
\def\@oddhead{{\small\lhead\ifnum\count0=\startpage ISSN 1464-8997 (on line)
1464-8989 (printed) \hfill {\lnum\number\count0}\else\ifodd\count0
\def\\{ }\ifx\theshorttitle\relax \thetitle \else\theshorttitle\fi\hfill
{\lnum\number\count0}\else\def\\{ and }{\lnum\number\count0}
\hfill\ifx\theshortauthors\relax 
\theauthors\else\theshortauthors\fi\fi\fi}}\def\@evenhead{@oddhead}
\def\@oddfoot{\small\lfoot\ifnum\count0=\startpage Published \publishdate:\qua\copyright\ \gtp\hfill\else
\gtm, Volume \thevolumenumber\ (\thevolumeyear)\hfill\fi}
\def\@evenfoot{@oddfoot}
\makeatother
\fi

\let\maketitlepage\maketitlep
\let\makeshorttitle\maketitlepage
\let\maketitle\maketitlepage

\newwrite\gtoutfile
\long\gdef\makeheadfile{  %%% start of definition of \makeheadfile
{\def\\{, }\def\s{ }
\immediate\openout\gtoutfile head.xxx
\immediate\write\gtoutfile{To: math@arxiv.org}
\immediate\write\gtoutfile{Subject: put OR rep NNNNN:ppppp}
\immediate\write\gtoutfile{--text follows this line--}
\immediate\write\gtoutfile{Proxy-for: \ifx\theasciiauthors\relax
\theauthors\else\theasciiauthors\fi\s<\ifx\theasciiemail\relax\theemail\else\theasciiemail\fi>}
\immediate\write\gtoutfile{\noexpand\\}
\immediate\write\gtoutfile{Authors: \ifx\theasciiauthors\relax
\theauthors\else\theasciiauthors\fi}
{\def\\{ }\immediate\write\gtoutfile{Title: \ifx\theasciititle\relax
\thetitle\else\theasciititle\fi}}
\immediate\write\gtoutfile{Subj-class: GT or SG, GR etc}
\immediate\write\gtoutfile{MSC-class: \theprimaryclass\ifx\thesecondaryclass\relax\else, \thesecondaryclass\fi}
\immediate\write\gtoutfile{Journal-ref: Geom. Topol. Monogr. \thevolumenumber\s
(\thevolumeyear) \startpage-\finishpage}
\immediate\write\gtoutfile{Comments: Published by Geometry and Topology Monographs at}
\immediate\write\gtoutfile{\s\s\s  http://www.maths.warwick.ac.uk/gt/GTMon\thevolumenumber/paper\thepapernumber.abs.html}
\immediate\write\gtoutfile{\noexpand\\}
\immediate\write\gtoutfile{}
\ifx\theasciiabstract\relax
\immediate\write\gtoutfile{\theabstract}\else
\immediate\write\gtoutfile{\theasciiabstract}\fi
\immediate\write\gtoutfile{}
\immediate\write\gtoutfile{\noexpand\\}
\immediate\write\gtoutfile{}
\immediate\closeout\gtoutfile}}  %%% end of definition of \makeheadfile

\def\maketitlepage{\maketitlep\makeheadfile}
\let\makeshorttitle\maketitlepage
\let\maketitle\maketitlepage

\volumenumber{4}
\volumename{Invariants of knots and 3-manifolds (Kyoto 2001)}
\volumeyear{2002}
\papernumber{5}
\pagenumbers{55}{68}
\received{30 November 2001}
\revised{8 April 2002}
\accepted{22 July 2002}
\published{19 September 2002}

\usepackage{amssymb,amsmath,rlepsf} 

\let\relabela\adjustrelabel

\newtheorem{thm}{Theorem}[section]  
        
\newtheorem{cor}[thm]{Corollary}

\theoremstyle{definition}

\newtheorem*{rem}{Remark}            
\newtheorem*{ackn}{Acknowledgement}

\newcommand\id{\operatorname{id}}

\newcommand\End{\operatorname{End}}

\newcommand\half{{\frac12}}

\newcommand\modZ {{\mathbb{Z}}}
\newcommand\modQ {{\mathbb{Q}}}
\newcommand\R{{\mathbb{R}}}
\newcommand\modR {{\mathcal{R}}}

\newcommand\onto\twoheadrightarrow

\newcommand\ho{{\hat\otimes }}

\newcommand\modv {{\mathsf v}}
\newcommand\tr{{\operatorname{tr}}}
\newcommand\clo{\operatorname{cl}}
\newcommand\ttr{\tilde\tr}
\newcommand\xqh[1]{{\widehat{#1[q]}}}
\newcommand\Zqh{{\xqh{\modZ }}}
\newcommand\hR{{\hat\modR }}

\begin{document}

\title{On the quantum $sl_2$ invariants of knots\\and integral homology spheres}
\asciititle{On the quantum sl_2 invariants of knots and integral homology spheres}
\authors{Kazuo Habiro}                  
\address{Research Institute for Mathematical Sciences\\Kyoto University, Kyoto, 606-8502, Japan}

\shorttitle{On the quantum $sl_2$ invariants}
\email{habiro@kurims.kyoto-u.ac.jp}

\begin{abstract}  
  We will announce some results on the values of quantum $sl_2$
  invariants of knots and integral homology spheres.  Lawrence's
  universal $sl_2$ invariant of knots takes values in a fairly small
  subalgebra of the center of the $h$-adic version of the quantized
  enveloping algebra of $sl_2$.  This implies an integrality result on
  the colored Jones polynomials of a knot.  We define an invariant of
  integral homology spheres with values in a completion of the Laurent
  polynomial ring of one variable over the integers which specializes
  at roots of unity to the Witten-Reshetikhin-Turaev invariants.  The
  definition of our invariant provides a new definition of
  Witten-Reshetikhin-Turaev invariant of integral homology spheres.
\end{abstract}

\asciiabstract{  
  We will announce some results on the values of quantum sl_2
  invariants of knots and integral homology spheres.  Lawrence's
  universal sl_2 invariant of knots takes values in a fairly small
  subalgebra of the center of the h-adic version of the quantized
  enveloping algebra of sl_2.  This implies an integrality result on
  the colored Jones polynomials of a knot.  We define an invariant of
  integral homology spheres with values in a completion of the Laurent
  polynomial ring of one variable over the integers which specializes
  at roots of unity to the Witten-Reshetikhin-Turaev invariants.  The
  definition of our invariant provides a new definition of
  Witten-Reshetikhin-Turaev invariant of integral homology spheres.}


\maketitle

\section{Introduction}
The purpose of this note is to announce some new results on the values
of quantum $sl_2$ invariants of knots and integral homology spheres.
We give a fairly small subalgebra of the center of the quantized
enveloping algebra $U_h(sl_2)$ of $sl_2$ in which Lawrence's universal
$sl_2$ invariant of knots takes values (Theorem~\ref{thm:1}).  This
implies a formula for the colored Jones polynomials of a knot
(Theorem~\ref{thm:4}).  We define an invariant $I(M)$ of integral
homology spheres $M$ with values in a completion of the Laurent
polynomial ring $\modZ [q,q^{-1}]$ which specializes at roots of unity to
the $sl_2$ Witten-Reshetikhin-Turaev (WRT) invariants.
(Theorem~\ref{thm:3}).  The definition of $I(M)$ leads to a new
definition of WRT invariant of integral homology spheres since we do
not use the WRT invariant in defining $I(M)$.  The invariant $I(M)$ is
as strong as the totality of the WRT invariants at various roots of
unity, and also as the Ohtsuki series (Theorem~\ref{thm:10}).  The
proofs, details and some generalizations will appear in separate
papers \cite{Habiro:ToAppear}.

\begin{ackn}
%  This work was presented at the project ``Invariants of Knots and
  Part of this work was presented at the project ``Invariants of Knots and
  $3$-Manifolds'' held at Research Institute for Mathematical
  Sciences, Kyoto University in September, 2001.  I am grateful to
%  Tomotada Ohtsuki for giving me the opportunity to give talks there.
  Tomotada Ohtsuki for giving me the opportunities to give talks there.
  I also thank Thang~Le for numerous conversations and correspondence.
\end{ackn}

\section{Lawrence's universal $sl_2$-invariant of links}

Let us recall the definition of the {\em universal
$sl_2$-invariant} introduced by Law\-rence \cite{Lawrence:univinv}.  She
actually introduced the invariant for more general Lie algebra, but we
will consider only the $sl_2$ case.

\subsection{The algebra $U_h(sl_2)$}
Let $U_h=U_h(sl_2)$ denote the quantized enveloping algebra of the Lie
algebra $sl_2$, i.e., the $h$-adically complete $\modQ [[h]]$-algebra
topologically generated by $H$, $E$, $F$ with the relations
$$
HE=E(H+2),\ HF=F(H-2),\ EF-FE=\frac{\exp({hH/2})-\exp({-hH/2})}{\exp({h/2})-\exp({-h/2})}.
$$
It is useful to introduce the elements
$$
  v = \exp({h/2})
$$
and 
$$
  K = v^{H} = \exp({hH/2}).
$$
For each $n\in \modZ $, set
$$
  [n] = (v^n-v^{-n})/(v-v^{-1}) \in  \modZ [v,v^{-1}].
$$
and for $n\ge 0$ set
$$
  [n]! = [1][2]\dots [n].
$$
The algebra $U_h$ has the structure of topological Hopf algebra given
by
\begin{gather}
  \Delta (H) = H\otimes 1+1\otimes H, \quad \epsilon (H) = 0,\quad S(H) = -H,\\
  \Delta (E) = E\otimes K+1\otimes E, \quad \epsilon (E) = 0,\quad S(E) = -EK^{-1},\\
  \Delta (F) = F\otimes 1+K^{-1}\otimes F, \quad \epsilon (F) = 0,\quad S(F) = -KF.
\end{gather}
The Hopf algebra $U_h$ has a universal $R$-matrix
$$
  R = v^{\half H\otimes H} \sum_{n\ge 0}
  v^{n(n-1)/2}\frac{(v-v^{-1})^n}{[n]!}(E^n\otimes F^n),
$$
and a ribbon element
$$
  r = K^{-1} \sum S(R_{(2)})R_{(1)},
$$
where we write $R=\sum R_{(1)}\otimes R_{(2)}$.

\subsection{Definition of the universal $sl_2$ invariant}

\begin{figure}
  \cl{\epsfbox{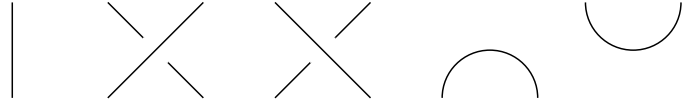}}
  \caption{The fundamental tangle diagrams}
  \label{fig:fundamental}
\end{figure}

Let $T$ be an oriented tangle diagram in $\R\times [0,1]$.  By small
isotopy we may assume that the critical points of the strings composed
with the ``height function'' onto $[0,1]$ are nondegenerate.  We may
also assume that on each crossing the two intersecting strings are not
critical.  As is well known, $T$ can be obtained from the {\em
  fundamental tangle diagrams} depicted in
Figure~\ref{fig:fundamental} by composition (i.e., pasting vertically)
and tensor product (i.e., pasting horizontally) up to isotopy of
$\R\times [0,1]$.

Let $T$ be a tangle diagram as above consisting of $l$ strings
$K_1,\dots ,K_l$ and $m$ circle components $K'_1,\dots ,K'_m$.  Choose a base
point $b_i$ of each $K'_i$ disjoint from the critical points and
crossings.

\begin{figure}
  \cl{\relabelbox\tiny
    \epsfxsize 4in \epsfbox{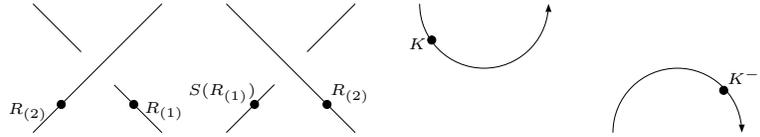}
\relabela <-1pt,0pt> {R1}{$R_{(2)}$}
\relabela <-2pt,0pt>  {R2}{$R_{(1)}$}
\relabela <-6pt,0pt>  {R3}{$R_{(2)}$}
\relabela <-1pt,0pt>  {S}{$S(R_{(1)})$}
\relabela <-1pt,0pt>  {K1}{$K$}
\relabela <-5pt,0pt>  {K2}{$K^{-1}$}
\endrelabelbox}
  \caption{How to put elements of $U_h$ on the strings}
  \label{fig:put}
\end{figure}

    We define below two elements $J'_T$ and $J_T$ by
    pretending for simplicity that we have $R=R_{(1)}\otimes R_{(2)}$ with
    $R_{(1)},R_{(2)}\in U_h$, though $R$ actually is an (infinite) sum
    of elements of the form $a\otimes b\in U_h\otimes U_h$.  The
    precise definition of $J'_T$ follows from 
    multilinearity.  We put elements of $U_h$ on the two strings near
    each crossing and on the string near each right-directed critical points as depicted in
    Figure~\ref{fig:put}.
  Then apply the
antipode $S$ to each element put on up-oriented strings.  For each
$K_i$, let $J_{(K_i)}$ be the expression obtained by reading from left
to right the elements on $K_i$ in the opposite orientation.
Similarly, let $J'_{(K'_i,b_i)}$ be the similarly obtained expression for
$K'_i$ reading from the base points $b_i$.  Then the expression
$$
  J'_T = \sum J_{(K_1)}\otimes \dots \otimes J_{(K_l)}\otimes J_{(K'_1,b_1)}\otimes \dots \otimes J_{(K'_l,b_l)}
$$
defines an element of $U_h^{\ho(l+m)}$, the $h$-adic completion of the
$(l+m)$-fold tensor product of $U_h$.

Let $I$ denote the $h$-adic closure of the $\modQ [[h]]$-submodule of
$U_h$ generated by the commutators $xy-yx$ for $x,y\in U_h$.  For $x\in U_h$,
let $[x]_I$ denote the coset $x+I\in U_h/I$.  For a tangle diagram $T$
as above, set $J_{(K'_i)}=[J_{(K'_i,b_i)}]_I$, which does not depend
on the choice of the base point of $K'_i$.  Then the expression
$$
  J_T = \sum J_{(K_1)}\otimes \dots \otimes J_{(K_l)}\otimes J_{(K'_1)}\otimes \dots \otimes J_{(K'_l)}
$$
defines an element of $U_h^{\ho l}\ho (U_h/I)^{\ho m}$, where $\ho$
denote the $h$-adically completed tensor product.  It is well known
that $J_T$ is invariant under isotopy of diagrams (fixing endpoints)
and framed Reidemeister moves, and defines an invariant of framed
tangles with total order on the set of components.  $J_T$ is called
the {\em universal $sl_2$ invariant} of $T$.

In the case of string knot $K$ (i.e., a $1$-string string link), the
universal $sl_2$ invariant $J_K$ of $K$ is contained in $U_h$.  It is
well known that $J_K$ is contained in the center $Z(U_h)$ of $U_h$.
It is also well known that $Z(U_h)$ is as a complete $\modQ [[h]]$-algebra
topologically freely generated by the element
$$
  c = FE + \frac{vK+v^{-1}K^{-1}-v-v^{-1}}{(v-v^{-1})^2}.
$$

\subsection{Integrality of the universal $sl_2$ invariant}
Let $C$ denote the well known central element of $U_h$
$$
  C = (v-v^{-1})^2 FE + vK + v^{-1}K^{-1} \in  Z(U_h)
$$
and set
$$
  \sigma _n = \prod_{i=1}^{n}(C^2-(v^{i}+v^{-i})^2)\in Z(U_h).
$$
We also set
$$
  q = v^2 \in \modQ [[h]],
$$
and regard $\modZ [q,q^{-1}]$ as a subring of $\modQ [[h]]$.

\begin{thm}
  \label{thm:1}
  Let $K$ be a string knot with $0$ framing.  Then there are unique
  elements $a_n(K)\in \modZ [q,q^{-1}]$ for $n\ge 0$ such that
  $$
  J_K = \sum_{n\ge 0} a_n(K) \sigma _n.
  $$
\end{thm}

For a knot $K$ with $0$-framing and an integer $n\ge0$, let
$a_n(K)\in\modZ[q,q^{-1}]$ denote the element determined by
Theorem~\ref{thm:1}.
If $K'$ is a string knot with $0$-framing with the closure 
 equivalent to $K$, then set $a_n(K') = a_n(K)$.

\section{Colored Jones polynomials of links}

\subsection{Finite dimensional representations of $U_h$}
By a {\em finite dimensional representation} of $U_h$ we mean as usual
a left $U_h$-module which is free of finite rank as a
$\modQ [[h]]$-module.  For each nonnegative integer $n$, there is exactly
one irreducible $(n+1)$-dimensional representation $V_{n+1}$ of $U_h$,
which corresponds to the unique $(n+1)$-dimensional representation of
$sl_2$.  The representation $V_{n+1}$ is defined as follows.  As a
$\modQ [[h]]$-module, $V_{n+1}$ is freely generated by the elements
$\modv _0,\modv _1,\dots ,\modv _n$.  The left action is given by
$$
   H\modv _i = (n-2i) \modv _i,\quad 
   E\modv _i = [n+1-i] \modv _{i-1},\quad F\modv _i = [i+1] \modv _{i+1},
$$
for $i=0,\dots ,n$, where we set $\modv _{-1}=\modv _{n+1}=0$.

For a finite dimensional representation $V$ of $U_h$ and a
$\modQ [[h]]$-module endomorphism $g\colon V\rightarrow V$, the {\em (left) quantum
trace} $\tr_q^V(f)$ of $f$ is defined by
$$
  \tr_q^V(f) = \tr(\rho _V(K)f),
$$
where $\rho _V\colon U_h\rightarrow \End V$ is the action of $U_h$ on $V$.  By abuse of
notation, we set for $x\in U_h$,
$$
  \tr_q^V(x) = \tr_q^V(\rho _V(x)),
$$
and call it the (left) quantum trace of $x$ in $V$.  We have
$$
  \tr_q^V(x) = \tr(\rho _V(Kx)).
$$
If $z\in Z(U_h)$, then $z$ acts on each $V_{n+1}$ as a scalar (i.e., an
element of $\modQ [[h]]$).  We have for any $v\in V_{n+1}$
$$
  z v = \frac{\tr_q^{V_{n+1}}(z)}{[n+1]} v.
$$
Here we have $\tr_q^{V_{n+1}}(1)=[n+1]$.  We have
$$
  \frac{\tr_q^{V_{n+1}}(z)}{[n+1]} = s_n(\varphi(z)),
$$
where $\varphi\colon U_h\rightarrow \modQ [H][[h]]$ is the $h$-adically continuous
$\modQ [[h]]$-linear map defined on the topological basis of $U_h$ by
$$
  \varphi(F^iH^jE^k) = \delta _{i,0}\delta _{k,0}H^j
$$
for $i,j,k\ge 0$, and $s_n \colon  \modQ [H][[h]]\rightarrow \modQ [[h]]$ denote the
$\modQ [[h]]$-algebra homomorphism defined by $s_n(g(H))=g(n)$.  The
restriction of $\varphi$ onto the $h$-adic closure $(U_h)_0$ of the
$\modQ [[h]]$-subalgebra spanned by $F^iH^jE^i$, $i,j\ge 0$, is a
$\modQ [[h]]$-algebra homomorphism known as the Harish-Chandra
homomorphism. It is well known that $\varphi$ maps the center
$Z(U_h)\subset (U_h)_0$ bijectively onto the $\modQ [[h]]$-subalgebra of
$\modQ [H][[h]]$ topologically generated by $(H+1)^2$.

\subsection{Colored Jones polynomial of links}
Let $L=(L_1,\dots ,L_l)$ be an ordered oriented framed link in $S^3$
consisting of $l$ components $L_1,\dots ,L_l$.  For nonnegative integers
$n_1,\dots ,n_l$, we can define the {\em colored Jones polynomial}
$J_L(V_{n_1+1},\dots ,V_{n_l+1})$ of $L$ associated with the ``colors''
$(n_1+1,\dots ,n_l+1)$ by
\begin{equation}
  \label{eq:3}
  J_L(V_{n_1+1},\dots ,V_{n_l+1}) = (\ttr^{V_{n_1+1}}\otimes \dots \otimes \ttr^{V_{n_l+1}})(J_L).
\end{equation}
Here $\ttr^{V_{n_i+1}}\colon U_h/I\rightarrow \modQ [[h]]$ is defined by
$$
\ttr^{V_{n_i+1}}([x]_I) = \tr(\rho _{V_{n_i+1}}(x)).
$$
(Usual definition of colored Jones polynomial involves braiding
operators on finite dimensional representations.  Our definition here is
equivalent to the usual one.)

We choose an $l$-component string link $T=(T_1,\dots ,T_l)$ such that the
closure of $T$ is ambient isotopic to $L$.  Then we have
$$
  J_L(V_{n_1+1},\dots ,V_{n_l+1}) =
  (\tr_q^{V_{n_1}+1}\otimes \dots \otimes \tr_q^{V_{n_l}+1})(J_T).
$$

\subsection{The case of knots of framing $0$}
Let $K$ be a string knot with $0$ framing.  Since $J_K\in Z(U_h)$,
$J_K$ acts on each representation $V_{n+1}$, $n\ge0$, as a scalar,
which we will denote by $J_K(V_{n+1})$.  It is well known that
$J_K(V_{n+1})\in\modZ[q,q^{-1}]$.  We have
$$
  J_{\clo(K)}(V_{n+1})=\tr_q^{V_{n+1}}(J_K)=[n+1]J_K(V_{n+1})
$$
and
$$
J_K(V_{n+1}) = s_n(\varphi(J_K)).
$$

\begin{thm}
  \label{thm:4}
  Let $K$ be a string knot with $0$ framing and let $n\ge 0$ be an
  integer.  Then we have
  \begin{equation}
    \label{eq:1}
    J_K(V_{n+1}) = \sum_{i=0}^n a_i(K) 
    \prod_{n+1-i\le j\le n+1+i,\ j\neq n+1}(v^j-v^{-j}).
  \end{equation}
  Note that this sum may be regarded as the infinite sum
  $\sum_{i=0}^\infty $ since the terms for $i>n$ vanishes.
\end{thm}

Theorem~\ref{thm:4} provides a new proof for Rozansky's integral
version \cite{Rozansky:MM} of the Melvin-Morton
expansion \cite{Melvin-Morton} of the colored Jones polynomials of
knots.  (We do {\em not} mean here that Theorem~\ref{thm:4} implies the
Melvin-Morton conjecture, proved in \cite{BarNatan-Garoufalidis:MM},
involving the Alexander polynomial.)  It follows from~\eqref{eq:1}
that
$$
  J_K(V_{n+1})
  = \sum_{i=0}^\infty  a_i(K)
  \prod_{j=1}^i (\alpha ^2 - (v^{j}-v^{-j})^2),
$$
where $\alpha =v^{n+1}-v^{-n-1}$.  The right hand side may be regarded as
an element of the completion ring
$$
  \varprojlim_{n} \modZ [q,q^{-1}, \alpha ^2] / (\prod_{j=1}^n(\alpha ^2 - q^j-q^{-j}+2)),
$$
with $\alpha ^2$ being regarded as an indeterminate.  There is a natural
injective homomorphism from this ring to the formal power series ring
$\modZ [[q-1, \alpha ^2]]$.

By expanding in powers of
$\alpha ^2$, we have
$$
  J_K(V_{n+1})
=\sum_{k=0}^\infty  \alpha ^{2k}(\sum_{i=k}^\infty  (-1)^{k-i} \tau _{i,i-k}a_i(K)),
$$
where 
$$
  \tau _{i,k} = \sum_{1\le p_1<\dots <p_k\le i}\prod_{r=1}^k(v^{p_r}-v^{-p_r})^2.
$$
It is not difficult to see that for each $k\ge 0$, the coefficient
$$\sum_{i=k}^\infty  (-1)^{k-i} \tau _{i,i-k}a_i(K)$$ of $\alpha ^{2k}$ defines an
element of the completion ring
$$
\Zqh=\varprojlim_{i}\modZ [q]/((q-1)(q^2-1)\cdots(q^i-1)).
$$
In particular, the constant term
$$
\sum_{i=0}^\infty  (-1)^{i} (\prod_{p=1}^i(v^p-v^{-p})^2)a_i(K)
$$
specializes to the Kashaev
invariants \cite{Kashaev,Murakami-Murakami} of $K$ by
substituting roots of unity for $q$.

\subsection{Examples}
Let $3_1^+$ (resp. $3_1^-$) denote the trefoil knot with positive
(resp. negative) signature, and let $4_1$ denote the figure eight
knot.
Then we have for each $n\ge 0$
\begin{gather}
  a_n(\text{unknot})  = \delta _{n,0},\\
  a_n(3_1^+) = (-1)^n q^{\half n(n+3)},\\
  a_n(3_1^-) = (-1)^n q^{-\half n(n+3)},\\
\label{eq:2}  a_n(4_1) = 1.
\end{gather}
A formula
for $4_1$ in \cite{Le:00} follows from \eqref{eq:2} and \eqref{eq:1}.

\subsection{The algebra $\modR $ and the basis $P'_n$}

% If $m,n\ge 0$, then we have a direct sum decomposion of left
If $m,n\ge 0$, then we have a direct sum decomposition of left
$U_h$-modules
$$
V_{m+1}\otimes V_{n+1} \cong \bigoplus_{|m-n|\le i\le m+n,\ i\equiv m+n\mod2}
V_{i+1}.
$$
The Grothendieck ring $\modR _\modZ $ of finite dimensional representations of
$U_h$ is freely spanned over $\modZ $ by $V_1=1,V_2,V_3,\dots $, and is
isomorphic to $\modZ [V_2]$.  For a commutative ring with unit $k$, set
$$
\modR _{k} = \modR _\modZ  \otimes _{\modZ } k = k[V_2].
$$
For $n\ge 0$, set
$$
  P_n = \prod_{i=0}^{n-1} (V_2 - v^{2i+1}-v^{-2i-1}) \in \modR _{\modZ [v,v^{-1}]}.
$$
We will also use the following normalizations
$$
P_n' = (v-v^{-1})^{-n} ([n]!)^{-1} P_n\in \modR _{\modQ (v)},
$$
and
$$
P_n'' = (v-v^{-1})^{-2n} ([2n+1]!)^{-1} P_n\in \modR _{\modQ (v)}.
$$
Let $L=(L_1,\dots ,L_l)$ be a framed link of $l$ components in $S^3$.
Extend \eqref{eq:3} multilinearly to define
$$
J_L(x_1,\dots ,x_l) \in \modQ (v^{1/2})
$$
for $x_1,\dots ,x_l\in \modR _{\modQ (v)}$.  If $L$ is algebraically split (i.e.,
the linking numbers are all $0$) and with all framings $0$, then we
have
$$
J_L(x_1,\dots ,x_l) \in \modQ (v).
$$

\begin{thm}
  \label{thm:2}
  If $K$ is a string knot with $0$ framing and if $n\ge 0$, then
  we have
  $$
    a_n(K) = J_{\clo(K)}(P''_n).
  $$
  In particular we have
  $$
    J_{\clo(K)}(P''_n) \in \modZ [q,q^{-1}].
  $$
\end{thm}

Let $\modR$ denote the $\modZ[v,v^{-1}]$-subalgebra of
$\modR_{\modQ(v)}$ generated by the elements $P'_n$ for $n\ge1$.  As
an $\modZ[v,v^{-1}]$-module, $\modR$ is freely generated by the
$P'_n$, $n\ge0$.
For each $m\ge 0$, let $\modR _m$ denote the $\modZ [v,v^{-1}]$-submodule of
$\modR $ spanned by $P'_m, P'_{m+1}, \dots $, which turns out to be an ideal
in $\modR $.  The following theorem is a generalization of
Theorem~\ref{thm:2} to algebraically split framed links.
\begin{thm}
  \label{thm:5}
  Let $L=(L_1,\dots ,L_l)$ be an algebraically split framed link of $l$
  components in $S^3$ with all framings $0$.
  If $x_1,\dots ,x_l\in \modR $, then we have
  $$
  J_L(x_1,\dots ,x_l) \in  \modZ [v,v^{-1}]
  $$
  If one of the $x_i$ is contained in $\modR _m$, $m\ge 0$, then we have
  $$
  J_L(x_1,\dots ,x_l) \in  \frac{(v-v^{-1})^m[2m+1]!}{[m]!}\modZ [v,v^{-1}]
  $$
\end{thm}

It follows from the first half of Theorem~\ref{thm:5} that if
$L=(L_1,\dots ,L_l)$ is an algebraically split framed link with
$0$-framings in $S^3$, then the $\modQ (v)$-multilinear map
$J_L\colon \modR _{\modQ (v)}\times \dots \times \modR _{\modQ (v)}\rightarrow \modQ (v)$ restricts to the
$\modZ [v,v^{-1}]$-multilinear map
\begin{equation}
  \label{eq:4}
  J_L\colon  \modR \times \dots \times \modR \rightarrow \modZ [v,v^{-1}],
\end{equation}
which induce the $\modZ [v,v^{-1}]$-linear map
\begin{equation}
  \label{eq:5}
  J_L\colon  \modR \otimes _{\modZ [v,v^{-1}]}\dots \otimes _{\modZ [v,v^{-1}]}\modR \rightarrow \modZ [v,v^{-1}].
\end{equation}
Set 
$$
\hR = \varprojlim_m \modR /\modR _m,
$$
which is a commutative $\modZ [v,v^{-1}]$-algebra.  $\hat\modR $ consists of the
infinite sums $\sum_{m\ge 0} b_m P'_m$, where $b_m\in \modZ [v,v^{-1}]$ for
$m\ge 0$.  It follows from the second half of Theorem~\ref{thm:5} that
$J_L$ in \eqref{eq:4} induces a $\modZ [v,v^{-1}]$-linear map
\begin{equation}
  J_L\colon \hR^{\ho l} \rightarrow  \widehat{\modZ [v]},
\end{equation}
where $\hR^{\ho l}$ denote the completion of the $l$-fold tensor
product $\hR\otimes _{\modZ [v,v^{-1}]}\dots \otimes _{\modZ [v,v^{-1}]}\hR$ with respect to
the natural tensor product topology induced by the completion topology
of $\hR$, and $\widehat{\modZ [v]}=\Zqh\otimes _{\modZ [q,q^{-1}]}\modZ [v,v^{-1}]$.
(Here recall that we set $q=v^2$.)

\section{A universal $sl_2$ invariant of integral homology\nl spheres}

In this section we define an invariant $I(M)\in \Zqh$ of integral
homology spheres $M$, which we call the ``universal $sl_2$-invariant'' of $M$, since
$I(M)$ is ``universal'' over the $sl_2$ WRT invariants at various
roots of unity.

\begin{rem}
  Recall that Le \cite{Le:aarhus} defined an invariant of closed
  $3$-manifolds $M$ with values in a ``functional space'' such that
  the $sl_2$ WRT invariants of $M$ recovers from it via certain ring
  homomorphisms.  However, his ``functional space'' is rather large,
  and since it involves complex functions, it does not give any
  information on the value of the WRT invariant at each root of unity.
\end{rem}

\begin{rem}
  The use of the word ``universal'' here is with respect to the roots
  of unity, but the use in Lawrence's $sl_2$ universal link invariant
  is with respect to finite dimensional representations.  We can unify
  these two invariants into a ``universal $sl_2$ invariant'' $I(M,L)$
  of links $L$ in integral homology spheres $M$.  This generalization
  is an easy modification of the definition of $I(M)$, and the details
  will appear in \cite{Habiro:ToAppear}.
\end{rem}

\subsection{The definition of the invariant $I(M)$}
Set
$$
\omega  = \sum_{i\ge 0} v^{\half i(i+3)}P'_i\in \hR,
$$
which is invertible in the algebra $\hR$ with the inverse
$$
\omega ^{-1} = \sum_{i\ge 0} (-1)^i v^{-\half i(i+3)}P'_i.
$$
Let $M$ be an integral homology $3$-sphere.  It is well known that
there is an algebraically split framed link $L=(L_1,\dots ,L_l)$ ($l\ge 0$)
in $S^3$ with all framings $\pm 1$ such that the surgery $(S^3)_L$ on
$S^3$ along $L$ is orientation-preserving homeomorphic to $M$.  Set
\begin{equation}
  \label{eq:6}
  I(L) = J_{L_0}(\omega ^{-f_1},\dots , \omega ^{-f_l})
\end{equation}
where $L_0$ denotes the framed link obtained from $L$ by changing all
the framings into $0$, and, for $i=1,\dots ,l$, $f_i=\pm 1$ denotes the
framing of the component $L_i$.  We have
$$
I(L) \in  \Zqh.
$$

\begin{thm}
  \label{thm:3}
  There is a well-defined invariant $I(M)$ of integral homology
  spheres $M$ with values in $\Zqh$ such that if $L$ is a
  algebraically split framed link in $S^3$ with all framings $\pm 1$,
  then we have $I((S^3)_L)= I(L)$.
\end{thm}

Theorem~\ref{thm:3} follows from Theorems~\ref{thm:6} and~\ref{thm:7}
below.
\begin{thm}[Conjectured by Hoste \cite{Hoste:Casson}]
  \label{thm:6}
  Let $L$ and $L'$ be two algebraically split framed links in $S^3$
  with all the framings $\pm 1$.  Then $L$ and $L'$ define
  orientation-preserving homeomorphic results of surgeries $(S^3)_L$
  and $(S^3)_{L'}$ if and only if $L$ and $L'$ are related by a finite
  sequence of {\em Hoste moves}, i.e., the usual Fenn-Rourke moves
  (surgery on unknotted component of framing $\pm 1$)
  {\em through} algebraically split framed links with $\pm 1$ framings.
\end{thm}

\begin{thm}
  \label{thm:7}
  The invariant $I(L)$ of algebraically split framed links $L$ in
  $S^3$ with $\pm 1$ framings is invariant under Hoste moves.
\end{thm}

\subsection{Specializations to the WRT invariants at roots of unity} 

For each root of unity $\zeta $, there is a well-defined ring homomorphism
$$
(-)|_{q=\zeta }\colon \Zqh\rightarrow \modZ [\zeta ],\quad f(q)\mapsto f(q)|_{q=\zeta }=f(\zeta ).
$$
For an integral homology sphere $M$ and a root of unity $\zeta $, let
$\tau _\zeta (M)$ be the WRT invariant of $M$ at $\zeta $ normalized so that
$\tau _\zeta (S^3)=1$.  (For definition of $\tau _\zeta (M)$, see
\cite{Kirby-Melvin}, but $\tau _r(M)$ for $r\ge 3$ defined there
corresponds to $\tau _{\exp(2\pi i/r)}$.  For the other primitive $r$th roots of
unity $\zeta $,
$\tau _\zeta (M)$ is obtained from $\tau _r(M)$ by the
automorphism of $\modZ [\zeta ]$ which maps $\exp(2\pi i/r)$ to $\zeta $.  For
$\zeta =\pm 1$, set $\tau _{\pm 1}(M)=1$.)

\begin{thm}
  \label{thm:8}
  For an integral homology sphere $M$ and a root of unity $\zeta $, we
  have
  \begin{equation}
    \label{eq:9}
    I(M)|_{q=\zeta } = \tau _\zeta (M).
  \end{equation}
\end{thm}

The proof of Theorem~\ref{thm:3} does not involve the existence proofs
of the variations of the WRT invariant in the literature.  Hence the
Theorems~\ref{thm:3} provides a new definition of the WRT invariant of
integral homology spheres via~\eqref{eq:9}.

\subsection{Consequences}
In the rest of this paper, we list some consequences to
Theorems~\ref{thm:3} and~\ref{thm:8}.

The following was first proved by H.~Murakami
\cite{H.Murakami}  for the case $\zeta $ is a root of unity of odd
prime order, and conjectured by Lawrence \cite{Lawrence:integrality2} in the
general case.
\begin{cor}
  \label{thm:9}
%  For an integral homology sphere $M$ and a root of unity, we have
  For an integral homology sphere $M$ and a root $\zeta$ of unity, we have
  $$
  \tau _\zeta (M) \in \modZ [\zeta ].
  $$
\end{cor}

By Lawrence's conjecture \cite{Lawrence:integrality} proved by
Rozansky \cite{Rozansky:lawrence}, the Ohtsuki series \cite{Ohtsuki:series} of an
integral homology sphere $M$ can be characterized as the formal power
series $\tau (M) \in  \modZ [[q-1]]$ such that for each root of unity $\zeta $ of
odd prime power order we have
\begin{equation}
  \label{eq:7}
  \tau (M)|_{q=\zeta } = \tau _\zeta (M),
\end{equation}
where both sides are regarded as the elements of $\modZ _p[\zeta ]$ with $p$
the odd prime such that the order of $q$ is a power of $p$.  Here
$\modZ _p$ denotes the ring of $p$-adic integers. 

Theorem~\ref{thm:8} provides a new proof of the existence of $\tau (M)$,
and moreover the following version of Lawrence's $p$-adic convergence
conjecture for $p=2$.
\begin{thm}
  If $\zeta $ is a primitive $2^m$th root of unity ($m\ge 1$), then we have
  \begin{equation}
    \label{eq:8}
    \tau (M)|_{q=\zeta } = \tau _\zeta (M) \in \modZ _2[\zeta ].
  \end{equation}
\end{thm}

Let
$$
\iota _1\colon \Zqh\rightarrow \modZ [[q-1]]
$$
be the homomorphism induced by $\id_{\modZ [q]}$.

\begin{thm}
  If $M$ is an integral homology sphere, then we have
  \begin{equation}
    \label{eq:10}
%    \iota (I(M))=\tau (M).
    \iota_1 (I(M))=\tau (M).
  \end{equation}
\end{thm}

Since $\iota _1$ is injective, $I(M)$ is as strong as $\tau (M)$.  The
injectivity of $\iota _1$ is also independently proved by Vogel.  We also
have the following.
\begin{thm}
  \label{thm:10}
  Let $M$ and $M'$ be two integral homology spheres.  Then the
  following conditions are equivalent.
  \begin{enumerate}
  \item $I(M) = I(M')$,
  \item $\tau (M) = \tau (M')$,
  \item $\tau _\zeta (M)=\tau _\zeta (M')$ for all roots of unity $\zeta $,
  \item $\tau _\zeta (M)=\tau _\zeta (M')$ for infinitely many roots of unity $\zeta $
    of prime power order.
  \end{enumerate}
\end{thm}

\begin{rem}
  For each root of unity $\zeta $, there is a natural homomorphism
  $$
  \iota _\zeta \colon \Zqh\rightarrow \modZ [\zeta ][[q-\zeta ]].
  $$
  For an integral homology sphere $M$, we may think of
  $\iota _\zeta (M)\in \modZ [\zeta ][[q-\zeta ]]$ as an ``expansion of the WRT invariants
  of $M$ at $q=\zeta $''.  This is a generalization of the Ohtsuki series
  to the expansion at $\zeta $.
\end{rem}

\Addresses\recd
\end{document}